\documentclass[12pt]{article}
\usepackage{amssymb,amsmath,amsfonts,mathrsfs,upgreek,txfonts,phonetic} 
\usepackage{mathabx}

\usepackage{graphicx} 
\usepackage{bm}

\usepackage[colorlinks=true]{hyperref}
\usepackage{comment}

\renewcommand{\marginpar}[2][]{} %attivarlo per eliminare i commenti a margine

\renewcommand\comment[1]{{\iffalse #1 \fi}}

\newtheorem{theorem}{Theorem}[section]

\newtheorem{lemma}{Lemma}[section]
\newtheorem{proposition}{Proposition}[section]

\newtheorem{remark}{Remark}[section]
\newtheorem{notation}{Notation}[section]
\newcommand{\diam}{{\,\rm diam\,}}
\newcommand{\io}{{\infty}}

\newcommand{\real}{ {\mathbb R}   }
\newcommand{\torus}{ {\mathbb T}    }
\newcommand{\integer}{ {\mathbb Z}   }
\newcommand{\complex}{ {\mathbb C}   }
\newcommand{\bB}{ {\mathbb B}   }

\newcommand{\cB}{ {\mathcal B}   }

\newcommand{\nP}{{\mathcal P}}

\renewcommand{\Im}{\, {\rm Im}\,}

\newcommand\beq[1]{ \begin{equation}\label{#1} }
\newcommand{\eeq}{ \end{equation} }
\newcommand{\beqno}{ \[ }
\newcommand{\eeqno}{ \] }
\newcommand\beqa[1]{ \begin{eqnarray} \label{#1}}
\newcommand{\eeqa}{ \end{eqnarray} }
\newcommand{\beqano}{ \begin{eqnarray*} }
\newcommand{\eeqano}{ \end{eqnarray*} }

\newtheorem{definition}{Definition}[section]
\newcommand\dfn[1]{ \begin{definition}\label{#1} }
\newcommand\edfn{ \end{definition} }
\newcommand\notat[1]{ \begin{notation} \label{#1} 
 }
\newcommand\enotat{\end{notation}}

\newcommand\rem{\begin{remark} 
\rm 
}
\newcommand\erem{\end{remark} %\vglue0.7truecm\noindent
}

\newcommand\equ[1]{{\rm (\ref{#1})}}
\newcommand{\nl}{{\smallskip\noindent}}
\newcommand{\giu}{{\medskip\noindent}}

\newcommand{\noi}{{\noindent}}

\newcommand{\qedeq}{\hskip.5truecm
\vrule width 1.7truemm height 3.5truemm depth 0.truemm}
 \newcommand\casialt[3]{ \left\{  \begin{array}{ll}
 {#1} & \mbox{ {\rm if} ${#2}$} \\
 {#3} & \mbox{ {\rm otherwise}}
 \end{array} \right.}

\newcommand{\e}{\varepsilon}

\renewcommand{\a }{\alpha }

\newcommand{\s }{\sigma }
\newcommand{\ii }{{\rm i} }
\renewcommand{\d }{\delta }

\newcommand{\f }{\varphi}

\renewcommand{\L }{\Lambda }
\newcommand{\m }{\mu }

\renewcommand{\t }{\tau }
\renewcommand{\o }{\omega }

\newcommand{\Z}{\mathbb{Z}}

\newcommand{\ttM}{{\rm  M}}

\newcommand{\ttH}{{\rm H}}
\newcommand{\ttD}{{\rm D}}

\def\N{\mathbb N}

\def\T{\mathbb T}

\def\const{{\, \rm const\, }}
\def\cc{{\, \rm c\,}}
\def\dst{\displaystyle}
\def\bks{\, \backslash\, }
\def\meas{{\rm\, meas\, }}
\def\Tp{{T_K^\perp}}

\newcommand\proiezione{\, {\rm p}}

\newcommand\ham{H_\e}
\newcommand\sa{\theta} %\newcommand\sa{\theta}
\newcommand{\tetta}{\vartheta }

\newcommand{\KO}{K}
\newcommand\ks[1]{K_{\rm o}(#1)}

\newcommand{\bs}{{\bar s}}

\newcommand{\norma}{\thickvert\!\!\thickvert}
\usepackage{appendix}

\newcommand\gen{{\cal G}^2_1}

\def\genKO{{\cal G}^n_{1,\KO}}

\newcommand\noruno[1]{  |#1|_{{}_1} }

\renewcommand\ln{\log}

\newcommand\hol{{\bB}}
\newcommand\palla{{\bf B_{\rm 1}}}

\newcommand\Ak{{\rm A}_k}

\newcommand\Dk{D^{1,k}}

\newcommand\ku{|k|_{{}_1}}
\newcommand\tc{\big|\,}

\newcommand\rescaling{\lambda}

\title{\bf 
On the measure of KAM tori in two degrees of freedom
}

\begin{document}

\author{ 
\footnotesize L. Biasco  \& L. Chierchia
\\ \footnotesize Dipartimento di Matematica e Fisica
\\ \footnotesize Universit\`a degli Studi Roma Tre
\\ \footnotesize Largo San L. Murialdo 1 - 00146 Roma, Italy
\\ {\footnotesize biasco@mat.uniroma3.it, luigi@mat.uniroma3.it}
\\ 
}

\maketitle

\begin{abstract}

\noi
A conjecture of Arnold, Kozlov and Neishtadt on the exponentially small  measure of the ``non-torus'' set in analytic systems with two degrees of freedom is discussed. 
\end{abstract}

{\small \tableofcontents}

%\date{}

\renewcommand{\=}{\coloneqq}

\section{Introduction and main result}
In this paper we consider real--analytic, nearly--integrable mechanical systems with two--degrees of freedom, namely, Hamiltonian systems governed by a Hamiltonian, in action--angle variables,  of the form
\beq{HM}
H_\e(y,x)\=\frac12 |y|^2 + \e f(x)\= \frac{y_1^2+y_2^2}{2} + \e\ f(x_1,x_2)\,,
\eeq
with 
$$y=(y_1,y_2)\in\real^2\,,\quad x=(x_1,x_2)\in \torus^2\= \real^2/(2\pi \integer)^2\,,\qquad
f:\torus^2 \to \real$$ 
real--analytic, $\e$ a small non negative  parameter. The phase space 
$\real^2\times \torus^2$  is  endowed with the standard symplectic form $dy_1\wedge dx_1+ dy_2\wedge dx_2$
so that the Hamiltonian flow induced by $H_\e$, 
$$
\phi_{H_\e}^t:(y_0,x_0)\in\real^2\times\torus^2\mapsto \big(y(t),x(t)\big)\= \phi_{H_\e}^t(y_0,x_0)\in 
\real^2\times\torus^2\ ,
$$
is the solution of standard Hamiltonian equations
$$
\left\{\begin{array}{l} 
\dot y = -\partial_x H_\e= -\e f_x\\
\dot x= \partial_y H_\e= y + \e f_y
\end{array}\right. \,, \qquad (y(0),x(0))=(y_0,x_0)\,.
$$
Such equation are equivalent to the Lagrangian Newtonian equations on $\torus^2$ with potential $f$, i.e.\footnote{As standard, dot denotes the derivative with respect to  ``time'' $t$  and $\partial_y=(\partial_{y_1},\partial_{y_2})$ and $\partial_x=(\partial_{x_1},\partial_{x_2})$ denote the gradients with respect to the variables $y$ and $x$.
}, 
$$
\ddot x= -\e f_x(x)\ ,\qquad \left\{\begin{array}{l}  x(0)=x_0\\
\dot x(0)=y_0\end{array}\right.\,.
$$ 

\nl
For $\e=0$, the system is integrable, the action variables $y_1$ and $y_2$ are integrals of the motions, and all trajectories  are simply given by $y(t)=y_0$ and $x(t)=x_0+ \o t$ where the frequency $\o$ coincides with the constant value  $y_0$. In particular the 2--tori $\{y_0\}\times \torus^2$ are all left invariant by the Hamiltonian flow
and whenever the ratio of the frequencies is an irrational number, such tori are spanned densely by any orbit.

\nl
As well known, according to classical KAM theory ``most'' integrable tori $\{y_0\}\times \torus^2$ persist for small 
$\e$ undergoing a small deformation  and fill any  bounded region of the phase space up to a set of measure at most $\sqrt\e$ (as $\e\to 0$);
these tori -- which are sometimes called {\sl primary tori} -- are Lagrangian graphs over $\torus^2$ and the motion is analytically conjugated to a translation by a Diophantine frequency\footnote{``Diophantine'' means that there exists $\a,\t>0$ such that 
$|\o\cdot k|= |\o_1k_1+\o_2k_2|\ge \a/|k|^\t$ for any non vanishing integer vector $k$.}
$\o$ on $\torus^2$;
(see,  \cite{AKN} for general information). 

\nl
This bound on the measure of the complement of primary tori
is sharp as it follows immediately by considering the trivial example
\beq{Pi}
H_\e= \frac{y_1^2+y_2^2}{2} + \e \cos x_1\ ,
\eeq
which governs the mechanics of  a simple pendulum with small gravity coupled with a free rotator. Indeed, this is an integrable 
system having different topologies for $\e=0$ and $\e>0$, and for $\e>0$ the measure of primary tori in any region $\{|y_i|\le R\}\times \torus^2$ with $\sqrt\e<R/2$, is given by  $(4\pi R)^2 (1- \frac{4}{\pi R}\sqrt\e)$. 

\nl
Of course, if one takes into account {\sl all} invariant tori, i.e., primary and {\sl secondary} tori (namely, the invariant tori that arise by effect of the perturbation and that in this trivial example correspond to
the $(y_1,x_1)$--librational orbits of the pendulum with initial data inside the separatrix $\{\frac12 y_1^2+\e \cos x_1=\e\}$), one has that the phase space of this integrable system is filled  by invariant Lagrangian tori,  up to a set of measure zero.

\nl
For general systems one does not expect to have a full set of invariant tori, however, 
Arnold, Kozlov and Neishtadt,  in  Remark 6.17 of \cite{AKN},  write:

\giu
{\sl It is natural to expect that in a generic (analytic) system with two degrees of freedom and with frequencies that do not vanish simultaneously the total measure of the ``non--torus'' set corresponding to all the resonances is exponentially small. However, this has not been proved.}

\giu
Indeed, we can prove the following result.

\giu
For $s>0$,  denote
\beq{toroseduto}
\torus^2_s\=\{x=(x_1,x_2)\in\complex^2 \tc\ |\Im x_j|<s\}/(2\pi \integer^2)\, ,
\eeq
and let $\hol_s^2$ be the Banach space of
{\sl real--analytic} functions on $\torus^2_s$ having zero average and finite  $\ell^\io$--Fourier norm\footnote{In this paper $x\cdot y$ denotes the inner product $x_1y_1+x_2 y_2$, $|x|$  the Euclidean norm $\sqrt{x_1^2+x_2^2}$ and $|x|_{{}_1}$ the 1--norm $|x_1|+|x_2|$; $f_k$ denotes the Fourier coefficient of order $k$, i.e., $\dst (2\pi)^{-2} \int_{\torus^2}f(x)e^{-\ii k\cdot x}\, dx$.}:
\beq{hol}
\hol_s^2\=\Big\{f=\sum_{k\in\Z^2 \atop k\neq 0} f_k e^{\ii k\cdot x} \ \tc \ \|f\|_s
\=\sup_{k\in\Z^2 \atop k\neq 0} |f_k| e^{|k|_{{}_1}s}<\io\Big\}\ .
\eeq

\nl
{\bf Theorem A}
{\sl Let $s>0$. There exists a set $\displaystyle\nP_s\subseteq \hol_s^2$, containing an open and dense set, such that the following holds.\\
Fix $0<r<R$, let $D\=\{y\in\real^2 \tc  r\le |y|\le R\}$ and consider the mechanical Hamiltonian system with phase space $D\times \torus^2$ and Hamiltonian $H_\e$ as in \equ{HM} with potential $f$ belonging to $\nP_s$. Then, 
there exists $\e_0,a>0$ small enough such that, whenever $0<\e<\e_0$, the Liouville measure of the complementary of $\phi^t_{H_\e}$--invariant  tori in the phase region $D$ is smaller than $R^2 \exp(-\const/\e^a)$.
}

\rem\label{secca}
(i) Notice that in the mechanical case the frequencies $\o_i\=\partial_{y_i} H_0=y_i$ vanish simultaneously only at $y=0$: this accounts for the annular shape of the action domain $D$  considered in the above theorem.

\nl
(ii) The exponent $a$ is computed in \cite{BC*}, where a detailed proof of the Theorem~A will appear.

\nl
(iii) The exponentially smallness of the ``non torus set'' (i.e.,  of the complementary of $\phi^t_{H_\e}$--invariant  tori) in two degrees of freedom is due to the fact that, in regions where the frequencies do not vanish simultaneously  (the origin, in the mechanical case)
{\sl there do not appear double resonances} (compare  Lemma~\ref{nodoppie} below).
%, a fact that simplifies significantly the analysis and in particular improves the needed normal forms. 

\nl
(iv) In three or more degrees of freedom,  multiple resonances instead are unavoidable and the exponential bound is in general no more valid. What one can prove  is the following 

\nl
{\bf Theorem} (\cite{BClin,BC})
{\sl
Consider a real--analytic nearly--integrable mechanical system with potential $f$, namely,  
a Hamiltonian system with  real-analytic Hamiltonian  
$$H_\e(y,x)=\frac12 \sum_{i=1}^n y_i^2  +\e f(x)\ ,$$  
$(y,x)\in\real^n\times\torus^n$ being standard action--angle variables.
For ``general non--degenerate'' potentials $f$'s there exists $\e_0,a>0$ such that, if $0<\e<\e_0$, then the Liouville measure of the complementary of $\phi^t_{H_\e}$--invariant  tori is smaller than $\e|\log \e|^a$.}

\giu
The class of ``general non--degenerate'' potentials  is the natural extension to higher dimension of the class $\nP_s$ defined in Sect.~\ref{sec:Ps} below. Also this theorem is in agreement (up to the logarithmic correction) with a conjecture by Arnold, Kozlov and Neishtadt\footnote{\label{toronto}
\cite[Remark 6.18, p. 285]{AKN}: ``{\sl It is natural to expect that in a generic system with three or more degrees of freedom the measure of the ``non--torus'' set has order $\e$". 
%Indeed, the $O(\sqrt\e)$--neighbourhoods of two resonant surfaces intersect in a
%domain of measure $\sim \e$. In this domain, after the partial averaging taking into account the resonances under consideration, normalizing the deviations of the ``actions'' from the resonant values by the quantity $\sqrt\e$, normalizing time, and discarding the terms of higher order, we obtain a Hamiltonian of the form $1/2(Ap, p) + V (q_1, q_2)$, which does not involve a small parameter (see the definition of the quantity $p$ above). Generally speaking, for this Hamiltonian there is a set of measure $\sim 1$ that does not contain points of invariant tori. Returning to the original variables we obtain a ``non--torus'' set of measure $\sim \e$.
}}.
\erem

\nl
In the rest of the paper,  we shall define $\nP_s$ and sketch the proof of Theorem~A.

\section{The generic set $\nP_s$}\label{sec:Ps}
Fix once and for all $s>0$. 

\nl
In this section  we define the generic set of potentials $\nP_s$. 

\nl
Denote by $\gen$ the ``generators'' of one--dimensional maximal lattices in $\integer^2$, i.e., 
\beq{cippalippa}
\gen\=\{k=(k_1,k_2)\in\integer^2: k_1>0\  {\rm and} \ \ {\rm gcd} (k_1,k_2)=1\}\cup \{(0,1)\}\ .
\eeq
Then,  the list of one--dimensional maximal lattices in $\integer^2$  is given by the sets $\integer k$ with $k\in \gen$
(explaining the name given to $\gen$).

\nl
Given a function $f\in \hol_s^2$ and given $k\in\gen$, we can project $f$, in Fourier space, on the lattice generated by $k\in\gen$ obtaining a function of the ``angle'' $k_1x_1+k_2x_2$, as follows
$$\sum_{j\in\integer} f_{jk} e^{\ii jk\cdot x}\eqqcolon F^k(k\cdot x)$$
where $\sa\to  F^k(\sa)$ is a real--analytic function on $\torus$ defined by
\beq{F^k}
F^k(\sa)= \sum_{j\in\integer\atop j\neq 0} f_{jk} e^{\ii j\sa}\,.
\eeq
One can, then,  decompose (in a unique way) the potential $f$ as sum of ``one dimensional'' functions of the angles $x\cdot k$, as $k\in\gen$:
\beq{visciola}
f(x)=\sum_{k\in\integer^2\atop k\neq 0} f_k e^{\ii k\cdot x} = \sum_{k\in \gen } F^k (x\cdot k)\,. 
\eeq
The functions $F^k$ will play a fundamental role in the forthcoming analysis.

\begin{definition}\label{fragola}
Let $0<\d\leq 1$ and let   
\beq{miguel}
\ks{\d}\=\cc  \max\big\{ 1\ ,\ \frac{1}{s}\ ,\ \frac{1}{s}\ \log\frac{1}{s \, \d} \big\}	\,,
\eeq
where $\cc>1$ is
 a suitable universal constant.
Denote by $\nP_s(\d)$   the set of functions in $\hol^2_s$ 
%there exists $ 0<\b<1$  such that 
%for every $k\in\integer^n_*$
such that,  for all $k\in\gen$ with   $\ku> \ks{\d}$, one has:

\begin{itemize}
\item[\footnotesize{\bf (P1)}]
$\displaystyle{  |f_k|\geq \d \ku^{-2}\ e^{-\ku s}\,,}$
\end{itemize}
while, for all $k\in\gen$ with  $\ku\le \ks{\d}$, one has:

\begin{itemize}
\item[\footnotesize{\bf (P2)}]
$\displaystyle{\  
%\min_{  |k|\leq \ks{\d} }\ \ 
\min_{\sa\in \torus} \ \big( |\partial_\sa F^k(\sa)|+|\partial^2_\sa F^k(\sa)|\big) >0}$\,;

\item[\footnotesize{\bf (P3)}]   $F^k(\sa_1)\neq F^k(\sa_2)$
for every $0\leq \sa_1<\sa_2<2\pi$ such that
$\partial_\sa F^k(\sa_1)=\partial_\sa F^k(\sa_2)=0$\,. 
\end{itemize}
Then, $\displaystyle\nP_s\=\bigcup_{\d>0} \nP_s(\d)$.
\end{definition}

\rem
(i) It is easy to produce  functions in $\nP_s(\d)$. Consider, for example, the function 
\beq{esempietto}
f(x)\= 2\d \sum_{k\in \gen} 
|k|_{{}_1}^{-{2}} e^{-|k|_{{}_1}s}\,  \cos (k\cdot x)\ .
\eeq
Such function 
has Fourier coefficients
\beqno
 f_k=\casialt{\dst\d |k|_{{}_1}^{-{2}} e^{-|k|_{{}_1}s}, }{\pm k\in\gen}{0 ,}
\eeqno
and Fourier projections 
$$
 F^k(\sa)=\d |k|_{{}_1}^{-{2}} e^{-|k|_{{}_1}s} \cos \sa\ .
$$
As it is plain,  $f\in \nP_s(\d)$.

\nl
(ii) The functions in   $\nP_s$ are  general in several  ways. \\
For example, from Proposition~3.1 of \cite{BCtopa},
 it follows easily that: 
\begin{itemize}
\item[(a)] {\sl $\nP_s$ contains an open and dense set in $\hol_s^2$}.

\item[(b)] {\sl $\nP_s$ is a prevalent set}\footnote{Recall that a Borel set $P$ of a Banach space $X$ is called {\sl prevalent} if there exists a compactly supported probability measure $\nu$ on the Borellians of $X$  such that $\nu(x+P)=1$ for all $x\in X$; compare, e.g., 
\cite{HK}.
\label{nurzia}}.

\item[(c)] The (weighted) Fourier map 
$$
j:f\in  \hol_s^2 \to \big\{f_k e^{|k|_{{}_1}s}\big\}_{k\in \gen}\in \ell^\io(\gen)
$$
yields  a natural isomorphisms between functions 
in $\hol_s^2$ and  bounded sequences of complex numbers supported on $\gen$.
\\
Denote by $\palla$ the closed ball of radius one in $\hol_s^2$ and by $\cB$ the Borellians in $\palla$.
\\
On $\palla$ one can introduce a natural (product) {\sl probability measure}, as follows. 
Consider, first,  the probability measure given by the  normalised Le\-bes\-gue--product measure 
on the unit closed ball of $\ell^\io(\gen)$, namely,   the unique probability measure $\mu$ on the Borellians of $\{z\in \ell^\io(\gen)\tc  |z|_{{}_\infty}\le 1\}$ such that,
given   Lebesgue measurable  sets $\Ak $ in the unit complex disk $D_1\=\{w\in\complex:\ |w|\le 1\}$
with $\Ak \neq D_1$ only for finitely many $k$, one has
 $$
\mu \Big(\prod_{k\in \gen} \Ak \Big)=
\prod_{\{k\in \gen:\, \Ak \neq D_1\}} 
\frac{1}{\pi}{\rm meas}(\Ak )\,
$$
where ``meas'' denotes  the Lebesgue measure on the unit complex disk $D_1$. 
\\
Then,  {\sl the isometry $j$ 
 induces a  probability measure 
$\mu_s$ on the Borellians $\cB$ and one has that 
$$\nP_s\, \cap \, \palla\in \cB\ ,\qquad{\rm and}\qquad \m_s(\nP_s\cap \palla)=1\ .
$$
}

\item[(d)] Assumption {\footnotesize{\bf (P3)}} is made in order to simplify (the quite technical and intricate) proofs but it is possible to obtain the main result also without such assumption.
\\
Assumption  {\footnotesize{\bf (P2)}} was used in \cite{Nei1} (see also \cite{Nei2}).

\end{itemize}

\erem

\section{Sketch of the proof of Theorem A}
Let $f\in \nP_s$ (Definition~\ref{fragola}), i.e.,  $f\in \nP_s(\d)$ for some $\d>0$, which will  henceforth be fixed.

\nl
In what follows, we denote by $c$ various (possibly different) constants, which may depend upon $s$, $\d$, $r$ and $R$.

\subsection{Small divisors and geometry of resonances}

Let $\a>0$ and $\KO\in\N$: $\a$ will measure the small divisors appearing and $\KO$ will be a Fourier cut--off. Later on these parameters will be suitably chosen as functions of $\e$ (see \equ{mora} below). \
In terms of these two parameters we shall describe the geometry of resonances. \\
Define  
\begin{itemize}

\item[\tiny $\bullet$] 
$\dst 
D^0
\coloneq \{ y\in D \ | \ 
|y\cdot k|\geq \a\,,\ \ \forall k\in\gen\,, \ |k|_1\leq\KO  \}$;

\item[\tiny $\bullet$]
$\dst 
D^{1,k}\coloneq \{ y\in D \ | \ 
|y\cdot k|<\a \}$, for $k \in \gen$;

\item[\tiny $\bullet$]
$\dst
D^1\coloneq \bigcup_{k\in\gen,\, |k|_1\leq\KO}
D^{1,k}\,;
$

\item[\tiny $\bullet$] For  $k\in\real^2\bks\{0\}$, denote by $\pi_k: \real^2\to \langle k\rangle\= \{t k\ |\ t\in\real\}$ the orthogonal projection 
onto the 1--dimensional vector space containing $k$, i.e., 
$$\dst \pi_k y\= \frac{y\cdot k}{|k|^2}\ k\,,$$
and by $\pi_k^\perp$ the orthogonal projection onto $\langle k\rangle^\perp$, the vector space orthogonal to $k$. Notice, that since we are in two space dimensions, $\langle k\rangle^\perp$ is the one--dimensional vector space containing $(k_2,-k_1)$, so that:
\beq{ovvieta}
\pi_k^\perp y \= y-\frac{y\cdot k}{|k|^2}\ k = 
\frac{y_1k_2-y_2k_1}{|k|^2}\ (k_2,-k_1)\,.
\eeq

\end{itemize} 

\rem\label{risonate}
(i) Recall that for the model at hand, frequencies $\o=\partial_y H_0$ and actions $y$ coincide.

\nl
(ii)
In the language of \cite{poschel}, $D^0$ is a $(\a,\KO)$--completely non resonant set, while $D^{1,k}$ is an $\a$--neighbourhood of an exact  
resonance $y\cdot k=0$ with $k\in\gen$ and $\ku\le \KO$; compare also Appendix~\ref{formosa}.

\nl
(iii) Obviously, by the definitions given, it follows immediately that
$$
D=D^0\cup D^1\,.
$$

\nl
(iv) For general ``geometry of resonances'' in the context of nearly--integrable Hamiltonian systems, see, e.g., 
\cite{nek}, \cite{poschel} and, more recently, \cite{GCB}. For a geometry of resonances specific for two--frequencies systems,
see \cite{Dirac}, \cite{Arnold} and \cite{Nei1}.
\erem

\begin{lemma}\label{nodoppie}
Let $\alpha\leq r/32\KO$,  $k\in\gen$ with $|k|\le \KO$. Let, also, $\ell\in\Z^2\bks  k\Z$ with $|\ell|\le 8\KO$. Then, 
$$
|y\cdot \ell|\geq \frac{r}{4|k|}\,,\qquad \forall\ y\in D^{1,k}\ .
$$
\end{lemma}
\begin{proof}
By \equ{ovvieta} and the definition of $D^{1,k}$, 
 $$
  |\pi_k^\perp y |\geq |y|-\frac{|y	\cdot k|}{|k|}
  > r- \frac{\a}{|k|}\geq \frac{r}{2}\,,
  $$ 
 and, observing that  $k_2 \ell_1-k_1 \ell_2\in\Z\bks\{0\}$ (since $\ell\notin k\Z$), 
  $$
  |\pi_k^\perp \ell |
    =\frac{|k_2 \ell_1-k_1 \ell_2| }{|k|}\geq \frac{1}{|k|}\,.
  $$
  Thus, (using again that $\langle k\rangle^\perp$ is one--dimensional), 
  \beqano
  |y\cdot \ell|&=& |\pi_k^\perp  y \cdot \pi_k^\perp  \ell + \pi_k  y \cdot\ell|\ge 
  |\pi_k^\perp  y \cdot \pi_k^\perp   \ell|-  |\pi_k  y \cdot\ell| 
\\
&=&    |\pi_k^\perp  y| \  |\pi_k^\perp   \ell|-  |\pi_k  y| \ |\ell| \ge  \frac{r}{2|k|}  - \a \ \frac{|\ell|}{|k|}\ge \frac{r}{4|k|}\ .
\qedeq
\eeqano
\end{proof}

\nl
{\sl From now on we fix:} 
\beq{mora}
\a\=r/2\,,\qquad\qquad
\KO\=\e^{-a}\,,
\eeq
where $0<a<1/6$ will be chosen later small enough.

\subsection{Averaging and normal forms}
In this section we construct suitable normal forms in the sets $D^0$ and $D^{1,k}$. The main tool is 
Proposition~4.1 of \cite{BCtopa}, which, for convenience of the reader,  is reported in Appendix~\ref{formosa}.

\nl
To describe  the  normal forms, we need to introduce proper norms.
\\
Given a domain ${\rm D}\subset \real^2$ and ${\rm r}>0$, we denote by ${\rm D_r}$ the complex neighbourhood 
$$
{\rm D_r}\=\{y\in\complex^2\tc |y-y_0|<{\rm r}\ ,\ {\rm for\ some}\ y_0\in {\rm D}\}\,;
$$
for a real--analytic function  ${\rm f}:   \T_{\rm s}^n\to \complex $  or ${\rm f}: {\rm D_r}\times \T_{\rm s}^n\to \complex$, we let, respectively, 
\beq{crostata}
\|{\rm f}\|_{\rm s}=\sup_{j\in \integer^n}  |{\rm f}_j |e^{|j|_{{}_1} {\rm s}}\ ,\qquad 
\|{\rm f}\|_{\rm r,s}=\sup_{j\in \integer^n} \sup_{y\in \rm D_r} |{\rm f}_j(y)|e^{|j|_{{}_1} {\rm s}}\ ,
\eeq
where ${\rm f}_j, {\rm f}_j(y)$ denote   Fourier coefficients.

\nl
For a given sublattice $\L\subseteq \integer^2$, we denote by $\proiezione_\L$ the Fourier--projection on $\L$: 
$$
\proiezione_\L f\=\sum_{k\in\L}f_k e^{\ii k\cdot x}\,.
$$

\subsubsection{Normal form on the non--resonant set $D^0$}
 
Set
$$
r_0\=\a/2\KO\,.
$$
then 
$$
|y\cdot k|\geq \a/2\,,\qquad
\forall y\in D^0_{r_0}
\,,\qquad \forall 0<|k|\leq \KO\,.
$$
From Proposition~\ref{pesce} it follows that,
for $\e$ small enough, there exists a symplectic change of variables
\begin{equation}\label{trota}
\phi_{0}: D^0_{r_{0}/2}\times \T^2_{s(1-2/\KO)} \to 
D^0_{r_{0}} \times \T^2_s
\,,
\end{equation}
such that\footnote{$ f^{{\rm o}}$ corresponds to $ f_{\varstar\varstar}$ in Proposition~\ref{pesce}. } 
%\footnote{To avoid the introduction of too many symbols, we denote again the new variables by $(y,x)$.} 
\begin{equation}\label{prurito}
\ham\circ\Psi_{0}
=\frac{|y|^2}2+\e g^{{\rm o}}(y) +
\e f^{{\rm o}}(y,x)\ ,  \qquad
\quad
\langle f^{{\rm o}}\rangle=0\,,
\end{equation}
where $\langle \cdot \rangle=\proiezione_{\{0\}}(\cdot)$ denotes the 
average with respect to the angles $x$, and: 
\beq{552}
\sup_{D^0_{r_{0}/2}}| g^{{\rm o}}-\langle f \rangle|
\leq
c\frac{\e \KO^2}{\a^2} \,, \qquad
%\sup_{D^0_{r_{0}/2}\times \T^2_{s(1-2/\KO)/2}}
\| f^{{\rm o}} \|_{r_{0}/2,s(1-2/\KO)/2}
\leq
e^{-\KO s/3}
\,.
\eeq

%\bigskip
%PARENTESI DA SPOSTARE
%
% Given $f(x)=\sum_{h\in\Z^2} f_h e^{\ii h x}$ and $k\in \gen$  we define the following {\bf 1d--Fourier projector} on $\T^1$
%\beq{pale}
% \pi_{k\integer} f \=F^k
%  \quad {\rm where}\quad
%F^k(\sa)\=\sum_{j\in \Z} f_{jk}e^{\ii j \sa}\ ,\qquad
%\sa\in\T^1
%\eeq
%$f_{jk}$ being the Fourier coefficient of $f$ with Fourier index $jk\in\integer^2$.
%
%
%
%\nl
%It is immediate to see that:
%
%\beq{dec}
%f(x)= \sum_{k\in \gen} F^k(k\cdot x)
%\eeq 

\bigskip

\subsubsection{Normal forms on simply--non--resonant sets $\Dk$}\label{panza}
Fix $k\in \genKO$ and 
let
\beq{rk}
r_k\= \frac{r}{32|k|\KO}
\eeq
then
$$
y\in \Dk_{r_k} \ , \quad \ell\in\Z^2\bks  k\Z\,,\ \ \
|\ell|\leq 8 \KO
\quad
\Longrightarrow
\quad
|y\cdot \ell|\geq \frac{r}{4|k|}\,.
$$
By Proposition~\ref{pesce}, with $(\a,K)$ replaced by $(\frac{r}{4|k|}, 8 \KO)$, 
we see that, for $\e$ small enough, there exists a symplectic change of variables
\begin{equation}\label{canarino}
\Psi_k: 
\Dk_{r_k/2}
\times \T^n_{s_\star } \to 
\Dk_{r_k} \times \T^n_s 
\,,\qquad s_\star\=s(1-1/\KO)
\end{equation}
such that\footnote{$ f^k$ corresponds to $ f_{\varstar\varstar}$ in Proposition~\ref{pesce}. } 
\beqa{Hk}
\ham\circ\Psi_k
&\eqqcolon& \frac{|y|^2}2+ \e G^k_0(y)  +
\e G^k(y,k\cdot x) +
\e f^k (y,x)
\eeqa
where  
\beq{sassone}
\langle G^k(y,\cdot)\rangle=0 \ ,\qquad  
\proiezione_{{}_{k\Z}} f^k=0 \ ,
\eeq
and\footnote{Beware that $F^k$ and $G^k$ are functions of one angle variable, while $f^k$  depends on two angle variables.}
\begin{equation}\label{552bis}
%\sup_{\Dk_{r_k/2}\times \T^1_{s(1-1/\KO)\noruno{k}}}
\sup_{\Dk_{r_k/2}} |G^k_0(y)|\ ,\ \|  G^k-F^k\|_{r_k/2, s_\star\ku} 
\leq
 %\tetta_k\=
 c \, {\e \ku^2 \KO^2}\,,
\qquad
\| f^k \|_{r_k/2,s_\star/2}
\leq
2 e^{-(4\KO-1) s}\ .
\end{equation}

\rem
The function $G^k(y,\theta)$ will be called the {\bf effective potential} since, disregarding the small remainder 
$f^k$, it governs the (integrable) Hamiltonian evolution at simple resonances.
\erem

\subsection{Exponential density of primary tori in $D^0\times \T^2$}\label{ventilatore}
In this brief section we  show how the exponential density of primary tori in the region  $D^0\times \T^2$ is an immediate consequence of the KAM Theorem, if one chooses suitably the parameter $\KO$ as a function of $\e$.

\nl
Indeed, 
we can apply the KAM Theorem~\ref{KAM} to the Hamiltonian in \equ{prurito} with $h(y)=|y|^2/2 + \e g^{{\rm o}}(y)$: in this case $h_{yy}={\mathbb I}+O(\e)$ and the perturbation 
$\e f^{{\rm o}}$ has norm bounded by (see \equ{552}) $\e e^{-\KO s/3}$. Therefore, recalling \equ{mora}, where we chose $\KO=1/\e^a$, one sees that the KAM condition
\equ{enza} is met for $\e$ small enough and that, by \equ{gusuppo},  the relative measure of Diophantine primary tori in $D^0\times \torus^2$ is at least 
\beq{mario}
1- \exp\Big({-\frac{s}{6\, \e^a}}\Big)\,.
\eeq

\subsection{The typical effective potential at simple resonances}
In the neighbourhoods $D^{1,k}$ of simple resonances, after the averaging of \S~\ref{panza}, the strategy is to put the integrable Hamiltonian\footnote{Integrable, since it depends only on the angle $Q=k\cdot x\in \torus^1$.}
$$
h=\frac{|y|^2}2+ \e G^k_0(y)  +
\e G^k(y,k\cdot x) 
$$
into action--angle variables, to check Kolmogorov's non--degeneracy
and then to apply  the KAM Theorem~\ref{KAM}. \\
To do this one has, first,  to understand the topological structure associated to the effective potentials $G^k$ for $\ku\le \KO$. 

\rem
(i) In the case\footnote{Recall \equ{miguel} for the definition of $\ks{\d}$.}  $\ku\le \ks{\d}$, the topology of the phase space of the effective integrable Hamiltonian can be quite arbitrary, as long as it is  non--degenerate, in the sense that 
the critical points of $\sa\mapsto G^k(y,\sa)$ are non--degenerate and at different energy levels (compare  {\footnotesize $\bf (P2)$},  {\footnotesize $\bf (P3)$} above). 

\nl
(ii) On the other hand, {\sl for $\ku>\ks{\d}$, 
all effective potentials $G^k$ have the same topological features of a pendulum}, as we shall briefly describe. \\
We stress that, while the case in (i) concerns a fixed (i.e., $\e$--independent) number of modes, the case $\ks{\d}<\ku\le K$ concerns a number of modes, which goes to infinity when $\e$ goes to zero. It is therefore essential to have {\sl unform control} of the case $\ks{\d}<\ku\le K$.

\nl
(iii)
{\sl From now on, to simplify the exposition, we shall consider only the case of simple resonances with $\ks{\d} <\ku\le \KO$. }
\\
The case $0<\ku\le \ks{\d}$, is similar but more complicated and we omit the details in the present sketch of proof.  
\erem

\nl
Thus, {\sl from now on, we fix $k\in \gen$ with  $\ks{\d} <\ku\le \KO$}.

\subsubsection{Uniform pendulum--like structure of the effective Hamiltonian ($\ku>\ks{\d}$)}
Because of the fast decay of Fourier modes due to analyticity,   $F^k$  (recall \equ{F^k}) has the form:
$$
F^k(\sa)=\big(f_k e^{\ii \sa }+ f_{-k}e^{-\ii\sa}\big)+ O(e^{-2\ku s})=
2|f_k|\, \cos (\sa+\sa_k)+ O(e^{-2\ku s})\,,
$$
for a suitable $\sa_k\in[0,2\pi)$. Recalling {\footnotesize $\bf (P1)$}, we can factor $|f_k|$, getting
$$
F^k(\sa)=
2|f_k|\Big( \cos (\sa+\sa_k)+ O(\ku^2 e^{-\ku s})\Big)\,.
$$
In fact,  these identities hold in a strong norm (e.g.,  in $\|\cdot\|_b$ with $b>1$; compare \equ{crostata}).  

\nl
Then, by \equ{552bis} and  {\footnotesize $\bf (P1)$},  one has\footnote{Notice that:  {\sl if  $\ 0<s'<s$ and $\langle f \rangle=0$, then $\|f\|_{s'}\le e^{s'-s} \|f\|_s$.}}:
$$
\frac1{|f_k|}\, \|G^k-F^k\|_{r_k/2,2} \le  c\, \frac1{|f_k|}  \, e^{-s\ku} \|G^k-F^k\|_{r_k/2,s_\star\ku}
\le c\, \ku^4 \KO^2 \e\,.
$$
Hence, recalling \equ{Hk}, and using again  {\footnotesize $\bf (P1)$}, one gets
\beq{cornogrande}
\ham\circ\Psi_k
=:
\frac{|y|^2}2+ \e G^k_0(y)+
2|f_k|\e
\Big( \cos(k\cdot x +\sa^{(k)})+
{\bf G}^k(y,k\cdot x)+
{\bf f}^k (y,x)
 \Big)
\eeq
with
\beq{themis}
\|{\bf G}^k\|_{r_k/2,2}\le c\,  \KO^6 \e=:\eta\ ,\qquad 
\|{\bf f}^k\|_{r_k/2,s_\star/2} \, \le e^{-5\KO s/2}\ .
\eeq
Recalling that in \equ{mora} we assumed  $a<1/6$, we get
\beq{alfio}
\eta=O(\e^{1-6a})\ll1\,.
\eeq

\subsubsection{Rescaling}
For the upcoming analysis it is convenient to make the rescaling\footnote{In the following, for ease of notation,  we shall sometimes  drop  the dependence on $k$, which has been fixed.}
\beq{rescaling}
y\to 
\rescaling
 y\,,\qquad {\rm where}\ \ \ 
\rescaling\= \sqrt{2|f_k|\e} 
\eeq
followed by a time--rescaling obtained by 
dividing the Hamiltonian by $\rescaling^2=2|f_k|\e$, so as to obtain the Hamiltonian 
\beq{cornopiccolo}
{\bf H}_k
\coloneqq
h_k(y)+
\big( \cos(k\cdot x +\sa^{(k)})+
{\bf G}^k(\rescaling\,  y,k\cdot x)+
{\bf f}^k (\rescaling\,  y,x)
 \big)\ ,
\eeq
where
\beq{maiella}
h_k(y)\=\frac{|y|^2}2+ \frac{1}{2|f_k|}G^k_0(\rescaling\, y)\ .
\eeq
\subsubsection{The fast angle $Q_2=k\cdot x$}
By  Bezout's Lemma we can find ${\bar k}=({\bar k}_1,{\bar k}_2)\in\Z^2$ with $|{\bar k}|_\infty\leq |k|_\infty$ 
such that
$$
{\bar k}_1k_1-{\bar k}_2k_2=1\,.
$$
Let
$$
A\=\left(\begin{matrix}
{\bar k}_1 & {\bar k}_2
\cr 
k_1 & k_2 \cr
\end{matrix}\right)\,.
$$
Applying the canonical transformation 
$$
\Psi_A: (P,Q) \mapsto (y,x)\,,
\qquad
y\=A^T P\,,\qquad
x\=A^{-1} Q
$$
and noting that $k\cdot x=Q_2$
we get
\beq{caffe}
{\bf H}_k \circ\Psi_A
=
h_k(A^T P)+
\big( \cos(Q_2 +\sa^{(k)})+
{\bf G}^k( \rescaling\,A^T P,Q_2)+
{\bf f}^k ( \rescaling\,A^T P,A^{-1} Q)
 \big)\,.
\eeq
The aim of this transformation is that, now, the effective potential 
$$
 \cos(Q_2 +\sa^{(k)})+
{\bf G}^k( \rescaling\,A^T P,Q_2)
$$
depends only on one angle, i.e. $Q_2$.

\rem\label{vento}
The norms of $A$ and $A^{-1}$ is proportional to  $|k|_\io$, and therefore the angle analiticity domain becomes $\T^2_{{s}/{cK}}$.
\erem
\subsubsection{Decoupling the kinetic energy}
However, this has the unpleasant cost that  the main part of the quadratic part in $P$ (the ``kinetic energy'')
$\frac12 |A^T P|^2$ is no longer diagonal.
In order to diagonalise it 
one can consider  the symplectic map
\beq{pajata}
\Psi_U: (p,q)\mapsto (P,Q)\,,
\qquad
P\=U p\,,\qquad
Q\=(U^{-1})^T q\,,
\eeq
where
\beqno
U\=\left(\begin{matrix}
1 & 0
\cr 
-{\bar k}\cdot k |k|^{-2} & 1 \cr
\end{matrix}\right)\
\,,
\eeqno
Indeed, using such a map, since $A^T U=[\pi_k^\perp {\bar k},k]$,  one finds 
$$
\frac12|A^T U p|^2=
\frac12  |\pi_k^\perp {\bar k}|^2 p_1^2 +
\frac12 |k|^2 p_2^2\,.
$$
However,  {\sl $\Psi_U$ does not yield a diffeomorphism on $\torus^2$}
as, in general, $\frac{{\bar k}\cdot k}{|k|^2}\in\mathbb Q$
 is not integer and, therefore,
$$
Q_1=q_1+\frac{{\bar k}\cdot k}{|k|^2} q_2
$$
is not well defined for $q_2\in\torus^1.$
Nevertheless,  applying $\Psi_U$ to the 
``effective Hamiltonian''
$$
h_k(A^T P)+
\big( \cos(q_2 +\sa^{(k)})+
{\bf G}^k( \rescaling\, A^T P,q_2)
 \big)
 $$
we get
\beq{pescatore}
\frac12  |\pi_k^\perp {\bar k}|^2 p_1^2 +
\frac12 |k|^2 p_2^2+
 W(p)+
\big( \cos(q_2 +\sa^{(k)})+
V(p,q_2)
\big)\,,
\eeq
where
$$
W(p)\=\frac{1}{2|f_k|}G_0^k( \rescaling\, A^T U p)\,,\qquad
V(p,q_2)\={\bf G}^k( \rescaling\, A^T U p,q_2)
$$
satisfy
\beq{vongole}
\sup_{{\bf D}^k_{{\bf r}_k}}
\|\partial^2_p W\|\leq \eta\,,\qquad
\|V\|_{{\bf r}_k,2}\le \eta\,,
\eeq
with\footnote{The fact that we can choose 
$r_k/4\rescaling\,|k|$
as new analyticity radius follows by \equ{rescaling}
and  estimating  the
operatorial norms of the matrices $A$ and $U$ as
$\|A\|\leq 2 |k|$ and $\|U\|\leq 2.$}: 
\beq{maiala}
{\bf D}^k \=\frac1{\rescaling\,}\, U^{-1}(A^{-1})^T \Dk\,,\qquad
{\bf r}_k\= \frac{r_k}{4\rescaling\,\, |k|}\geq 1
\eeq
for $\e$ small enough (recall \equ{mora},
\equ{rk} and
 \equ{rescaling}).
 
 \subsubsection{Action--angle variables}
Since the ``effective Hamiltonian'' in \equ{pescatore}
does not depend on the angle $q_1 $, the action $p_1$
is an integral of motion and plays the role of a parameter.
Then, disregarding the dynamically irrelevant term  
$\frac12  |\pi_k^\perp h|^2 p_1^2$, we study the 
``pendulum-like Hamiltonian''
$$
H_{\rm pend}(p_2,q_2;p_1)
\=\frac12 |k|^2 p_2^2+
 W(p)+
\Big( \cos(q_2 +\sa^{(k)})+
V(p,q_2)
\Big)\,.
$$
$H_{\rm pend}$ is a one dimensional Hamiltonian depending on the
parameter $p_1$ and, therefore, it is integrable introducing
suitable action angle variable.

\nl
The separatrix divides the phase of $H_{\rm pend}$  into  three ($p_1$--dependent) open regions: 
${\mathcal D}_+$,  above the separatrix, 
${\mathcal D}_-$,  below the separatrix, and  
${\mathcal D}_0$,  inside the separatrix (excluding the elliptic equilibrium), which will contain the (projection of) the secondary tori, i.e.,  those Lagrangian tori, which are not graphs over the angles.
\\
Next, we construct, in each region, action--angle variables $(p_2,q_2)$ through $p_1$-dependent
symplectic transformations
\beq{vaccinara}
p_2=p_2^\s(I_2,\f_2;p_1)\,,\qquad
q_2=q_2^\s(I_2,\f_2;p_1)\,,
\eeq
with $\s=+,-$ or $0$, 
such that, in the new variable $H_{\rm pend}$ reads
\beq{gnocchi}
H_{\rm pend}(p_2^\s,q_2^\s;p_1)\eqqcolon E^\s(p_1,I_2)
\eeq
(which is integrable).
Note that the maps in \equ{vaccinara}
can be easily completed into symplectic transformations
$$
\Psi_{\rm aa}^\s:(I_1,I_2,\f_1,\f_2)\mapsto
(p_1,p_2,q_1,q_2)
$$
fixing $p_1=I_1$.

\nl
It is important to remark that, even though  $\Psi_U$ (defined in \equ{pajata})
is not well defined on the angles, 
the composition 
$$
\Psi_U\circ\Psi_{\rm aa}^\pm\circ \Psi_U^{-1}
$$
is instead well defined.

\nl
On the other hand, in the region ${\mathcal D}_0$, in view of the different topology, it is actually enough to consider the symplectic transformation
$$
\Psi_U\circ\Psi_{\rm aa}^0\,,$$
which is well defined.

\nl
%For simplicity, let us consider the inner region ${\mathcal D}_0$. 
In the variables $(I,\f)$ the Hamiltonian takes the form 
$$
{\bf h}_k(I)+{\bf f}_k(I,\f)\ ,
$$
where
\beq{gransasso}
{\bf h}_k(I)\=\frac12  |\pi_k^\perp {\bar k}|^2 I_1^2
+E^\s(I)\ ,\qquad {\bf f}_k =O\big(\exp(-5\KO s/2)\big)\,.
\eeq
\subsubsection{Kolmogorov's non--degeneracy}
In order to apply the KAM Theorem~\ref{KAM} to such Hamiltonian, we need to show that ${\bf h}_k$
twists, namely, that the determinant of its Hessian
is bounded away from zero.

\rem
Notice that, recalling \equ{vongole},  for $\eta= 0$, $E^\s(I)|_{\eta=0}$ reduces to the pendulum (in action variables) and the twist  can be checked by direct computations
As far as one stays away from the separatrix,  one can still check the twist  perturbatively. However, we need estimates in regions which  are {\sl exponentially (in $1/\e$) close to the separatrix} and this regime is no longer perturbative, as we are going to explain.
\\
Indeed, denoting by $z$ the distance in energy
from the separatrix, it can be shown that, asymptotically as $\eta,z\to 0$, one has (up to multiplicative $|\log z|^b$--corrections)
\beqa{citronella}
\det \partial_I^2 {\bf h}_k & \cong & \det
\left(\begin{array} {cc} 
|\pi_k^\perp {\bar k} |^2+ 
O(\eta/z)
& 
O(\eta/z) 
\\ 
O(\eta/z) 
& 
c_0/z  
\end{array}\right)\nonumber\\
&=&\frac{c_1}{z}+\frac{O(\eta)}{z^2}
\eeqa
with $c_1= |\pi_k^\perp {\bar k}|^2 c_0\neq0$,
and, since $z$ can be much smaller than $\eta$, we see that the evaluation in \equ{citronella} turns into a  
{\sl singular perturbation problem}, and hence cannot be handled by usual perturbation techniques.
\erem

\nl
To overcome this problem, we consider the inverse of the function $I_2\mapsto E=E^\s(I_1, I_2)$, parameterised by $I_1$: let us call it $I_2^\s(z;I_1)$, where $z\=E-E_0$, $E_0=E_0(I_1)$ being the energy of  the separatrix.

\nl
Now, one can prove that
$$
I_2^\s(z;I_1)=\phi^\s(z;I_1)+ \chi^\s(z;I_1) \ z \log z  \,,
$$
with $\phi^\s$ and $\chi^\s$ analytic in $z$ near the origin. 
\\
By using analyticity arguments, we can then show that:

\nl
{\sl  For any $\theta>0$ small enough,
up to a region $\theta$--bounded away from separatrices and of measure of order
$\theta^{c_1}$ for some $0<{c_1}<1$,
the following estimates hold 
uniformly in $|k|\le K$:
\beq{brega}
\|\partial^2_I {\bf h}_k \|\le \frac1{\theta}\ ,\qquad \qquad |\det \partial^2_I {\bf h}_k|\ge \theta\ .
\eeq
}

\nl

\subsection{Exponential density of primary and secondary tori in $D^1\times \T^2$}

In the  region where \equ{brega} holds,
 we can apply the KAM Theorem~\ref{KAM} with $d=\theta$, ${\rm M}=1/\theta$, 
$\mu=\theta^3$, $\e_0=O\big(\exp(-5Ks/2)\big)$, recall \equ{gransasso}, 
${\diam}{\rm D}\leq K/c\rescaling$,
 $r=\theta/c \KO$ and\footnote{Recall Remark \ref{vento}.} $s=1/c\KO$. 
Then 
$
\epsilon\leq e^{-5Ks/2}/\theta^2
$
and the KAM condition in \equ{enza}
is satisfied choosing 
\beq{theta}
\theta= \exp( - c_2 /\e^a)\,,
\eeq
for a suitable $c_2$ small enough and $\e$
small enough.
Since $C$ in \equ{ciccione} is bounded, for
$\e$
small enough,  by
$1/\rescaling^2 \theta^{22}\leq e^{Ks}/\theta^{23}$
(recall \equ{rescaling}),
then the measure of the complement of invariant tori is bounded, recalling \equ{gusuppo}, by
$$
C \exp(-5K s/4)
\leq \theta^{-24} \exp(-5K s/4)\le \exp(-K s/4)\le \theta\,,
$$  
for $\e$ small enough.
In conclusion, 
%in the region  inside separatrix 
the measure of the complement of invariant tori is bounded by $2\theta.$
%The same bound holds also in the upper and lower regions.
Recalling \equ{maiala}
we have that in the starting domain $\Dk\times\T^2$
  the measure of the complement of invariant tori is bounded by
 $$
\frac{1}{c} K^2 \rescaling^2\theta\,.
 $$
 Then,
 in the whole region $D^1\times\T^2$
 the measure of the complement of invariant tori is bounded by
 $$
\frac{1}{c} K^4 \rescaling^2\theta
\leq
\frac{1}{c'} K^4 \e \theta
\leq \theta
 $$
 for $\e $
 small enough.
 \\
 This last estimate, recalling the definition \equ{theta}, together with  the estimates of \S~\ref{ventilatore},  
 concludes the proof of Theorem A.

\appendix

\section{Appendix: Normal forms and  KAM}

\subsection{A normal form lemma}\label{formosa}

\nl
The following normal form lemma is proven in \cite[\bf Proposition~4.1]{BCtopa}. Before stating it we need some definitions.

\begin{itemize}
\item[\small $\bullet$] For  functions
 $f:D_r\times \torus^n_s\to\complex$ we set
\begin{equation}\label{canuto}
\norma f\norma_{D,r,s}=
\norma f\norma_{r,s}\=
\sup_{y\in D_r}\sum_{k\in\Z^n} 
|f_k(y)| e^{|k|_{{}_1}s}\,.
\end{equation}

\nl 
The norms $\|\cdot\|_{r,s}$ and  $\norma \cdot\norma_{r,s}$ are not equivalent, however the following relation holds

\beqa{battiato3}
\|f\|_{r,s}\leq  \norma f\norma_{r,s}&\le& 
(\coth^n(\s/2)-1)\|f\|_{r,s+\s}\\
&\leq&
(2n/\s)^n\|f\|_{r,s+\s}\,.\nonumber
\eeqa

\item[\small $\bullet$]
Given an integrable Hamiltonian $h(y)$, positive numbers $\a,K$ and a lattice $\L\subset \Z^n$, a (real or complex) domain $U$ is $(\a,K)$ non--resonant modulo $\L$ (with respect to $h$) if  
\begin{equation}\label{panzapiena}
|h'(y)\cdot k |\ge \a\ ,
\ \ \ \forall\ y\in U\ , \forall\  k\in \integer^n\bks \L\ , \    \noruno{k}\leq K\ .
\end{equation}

\item[\small $\bullet$] 
Given $f(y,x)=\sum_{k\in\Z^n}f_k(y)e^{\ii k\cdot x}$ and   a sublattice $\L$ of $\Z^n$, 
we denote by $\proiezione_\L$ the projection on the Fourier coefficients in $\L,$ namely
$$
\proiezione_\L f\=\sum_{k\in\L}f_k(y)e^{\ii k\cdot x}\,.
$$
and by $\proiezione_\L^\perp$ its ``orthogonal'' operator (projection on the Fourier modes in $\integer^n\bks\L$):
$$
\proiezione_\L^\perp f\=\sum_{k\notin\L}f_k(y)e^{\ii k\cdot x}\,.
$$

\end{itemize}

\begin{proposition}[\cite{BCtopa}]\label{pesce}\ \\
Let $r,s,\a>0$,  $K\in\N,$ $K\ge 2$, $D\subseteq \real^n$,  and let $\L$ be a lattice of $\integer^n$.
Let 
$$H(y,x)=h(y)+f(y,x)$$
 be real--analytic on  
$ D_r \times \T^n_s$ with $\norma f\norma _{r,s}<\infty.$
Assume that $D_r$ is
($\a$,$K$)--non--resonant modulo
%\footnote{In case $\L=\{0\}$, one also says that $D_r$ is completely $(\a,K)$ non--resonant.}
$\L$ 
 and that
$$
\tetta_\varstar \= \frac{2^{11} K^2}{\a r s}\, \norma f\norma _{r,s}  <1\,.
$$
Then,  there exists a real--analytic 
 symplectic change of variables
$$
\Psi: (y',x')\in D_{r_\varstar}\times \T^n_{s_\varstar}\  \mapsto\ 
(y,x)\in D_r \times \T^n_s \, \quad {\rm with} 
\quad
r_\varstar\=r/2\,,\ \ 
s_\varstar\=s(1-1/K)
$$
satisfying
$$
\noruno{y-y'}\leq \frac{\tetta_\varstar}{2^7 K} r\,,\qquad
\max_{1\leq i\leq n}|x_i-x'_i|
\leq 
\frac{\tetta_\varstar}{16 K^2} s\,,
$$
and such that 
\beq{senzanome}
H\circ\Psi=h+ f^\flat +
f_\varstar\,,   \qquad f^\flat \= \proiezione_\L f+ \Tp\proiezione_\L^\perp f
\eeq
with 
$$
%\norma \proiezione_\L f_\varstar\norma _{r_\varstar,s_\varstar}\,,\  \norma \Tp\proiezione_\L^\perp f_\varstar\norma _{r_\varstar,s_\varstar}
\norma f_\varstar\norma _{r_\varstar,s_\varstar}
\leq 
\frac{1}{K}\tetta_\varstar\norma f\norma_{r,s}
\,,\qquad 
\norma T_{K}\proiezione_\L^\perp f_\varstar\norma _{r_\varstar,s_\varstar} 
\leq
(\tetta_\varstar/8)^{{K}}\frac{8}{eK}\norma f\norma _{r,s}\,.
$$
Moreover,  re-writing \equ{senzanome} as 
$$
H\circ \Psi=h+g+ f_{\varstar\varstar}\qquad \mbox{where}\quad \proiezione_\L g=g\ ,\quad \proiezione_\L f_{\varstar\varstar}=0\,,
$$
one has
$$
\norma g- \proiezione_\L f\norma _{r_\varstar,s_\varstar}\le 
\frac{1}{K}\tetta_\varstar \norma f\norma _{r,s}\ ,\qquad
\norma f_{\varstar\varstar}\norma _{r_\varstar,s/2}
\le 
 2 e^{-(K-2)\bs}\norma f\norma _{r,s}\,,
$$
where
$$
\bs\=\min\left\{
\frac{s}{2},\,
\ln\frac{8}{\tetta_\varstar}
\right\}\,.
$$

\end{proposition}

\rem
The main point of Proposition \ref{pesce} concerns the analyticity domain in the angular variables of the renormalised Hamiltonian, which is close to optimal. Indeed, the Fourier
coefficients of the new Hamiltonian {\sl are shown to decay at the exact same exponential rate as the Fourier coefficients of the original Hamiltonian}, at least up to order $K$,  and this fact  plays a crucial role in our analysis.
\erem

\subsection{A KAM Theorem}

\begin{theorem}\label{KAM}
Let $r,s>0$, $n\ge 2$, $\ttD\subseteq \real^n$ be a bounded set and 
$H(y,x)=h(y)+f(y,x)$
 be a real--analytic Hamiltonian on  
$ \ttD_r \times \T^n_s$, such that
$$
\ttM\= \sup_{\ttD_r} |h_{pp}|<+\io\,,\qquad d\=  \inf_\ttD|\det h_{pp}| >0\,,\qquad
\e_0\=\sup_{\ttD_r \times \T^n_s} |f|<+\io \,.
$$
\begin{equation}\label{dupa}
\end{equation}
Let also 
\beq{bove}
\mu\=\frac{d}{\ttM^n}\,,
\eeq
and fix $\t> n-1$. \\
Then, there exists positive constants $c<1$  depending only on $n$ and $\t$ such that, if  
\beq{enza} 
\epsilon\= \frac{\e_0}{\ttM r^2}\le c \, \m^8\ s^{4\tau+8}\ , 
\eeq
then the following holds.
Define
\begin{equation}\label{nicaragua}
\a\=\frac{\ttM r}{ \m\, s^{3\tau +6}} \,\sqrt\epsilon\ ,\qquad
\hat r\= \m^2 r\ ,\qquad r_\epsilon\=\frac{1 }{c}\, \frac{\sqrt\epsilon\, r}{\m}\ .
\end{equation}
Then, there exists a positive measure 
%``non--torus set'' 
set ${\cal T}_\a\subseteq \ttD_{\hat r}\times \torus^n$ formed by ``primary''  Kolmogorov's tori; more precisely, 
for any point $(p,q)\in{\cal T}_\a$, $\phi^t_H(p,q)$ covers densely an $\ttH$--invariant, analytic, Lagrangian torus,  with $\ttH$--flow analytically conjugated to a linear flow with $(\a,\t)$--Diophantine frequencies 
$\o=h_p(p_0)$, for a suitable $p_0\in\ttD$; each of such tori is a graph over $\torus^n$ $r_\epsilon$--close
to the unperturbed  trivial graph  $\{(p,\theta)=(p_0,\theta)|\ \theta\in \torus^n\}$.\\
Finally, 
the Lebesgue outer measure of $(\ttD\times\torus^n)\bks {\cal T}_\a$ is bounded by:
\beq{gusuppo}
\meas \big((\ttD\times\torus^n)\bks {\cal T}_\a\big) \le  C\, \sqrt\epsilon
\eeq
with
\beq{ciccione}
C\=
\big(\max\big\{ \m^2 r\, ,\, \diam\, \ttD\big\}\big)^n \cdot \frac{1}{c \m^{n+5}\ s^{3\tau +6}}\ .
\eeq
\end{theorem}

\rem
(i) Theorem \ref{KAM} is an immediate consequence of 
Theorem~1 in \cite{BCKAM} (actually, it is just a slightly simplified version of it).

\nl
(ii) 
Notice that $\m\le 1$: in fact, since the eigenvalues of $h_{pp}$ are bounded in absolute value by $\|h_{pp}\|\le \ttM$,  one has that
$d\le \sup_\ttD |\det h_{pp}|\le \ttM^n$.

\nl
(iii) The main point of Theorem~\ref{KAM}  is to have a quantitative smallness condition with explicit dependence on the domain $\ttD$: This is  important  for our application, since  domains (after rescalings and changes of variables) may become very large.
\erem

\end{document}